\documentclass{article}

\usepackage{arxiv}

\usepackage[utf8]{inputenc}
\usepackage{graphicx}
\usepackage{xcolor}

\usepackage[colorlinks=true, allcolors=blue]{hyperref}

\usepackage{amsmath}
\usepackage{amsfonts}
\usepackage{bm}

\usepackage{booktabs}
\usepackage{multirow}
\usepackage{diagbox}

\title{Vibroacoustic simulations of acoustic damping materials using a fictitious domain approach}

\author{ \href{https://orcid.org/0000-0001-7015-8928}{\includegraphics[scale=0.06]{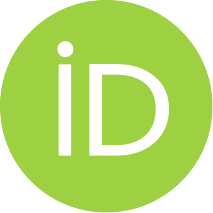}\hspace{1mm}Lars Radtke} \\
	Institute for Ship Structural Design and Analysis (M-10) \\
	Numerical Structural Analysis with Appl. in Ship Technology \\
	Hamburg University of Technology \\
	Hamburg, 21073 \\
	\texttt{lars.radtke@tuhh.de} \\
	\And
	\href{https://orcid.org/0000-0001-5937-6158}{\includegraphics[scale=0.06]{orcid.pdf}\hspace{1mm}Paul Marter} \\
	Institute~of~Mechanics \\
	Otto von Guericke University Magdeburg \\
	Magdeburg, 39106 \\
	\texttt{paul.marter@ovgu.de} \\
	\And
	\href{https://orcid.org/0000-0002-2396-4035}{\includegraphics[scale=0.06]{orcid.pdf}\hspace{1mm}Fabian Duvigneau} \\
	Institute~of~Mechanics \\
	Otto von Guericke University Magdeburg \\
	Magdeburg, 39106 \\
	\texttt{fabian.duvigneau@gmail.com} \\
		\And
	\href{https://orcid.org/0000-0001-8774-9732}{\includegraphics[scale=0.06]{orcid.pdf}\hspace{1mm}Sascha Eisentr\"ager} \\
	Institute~for~Mechanics \\
	Technical University of Darmstadt \\
	Darmstadt, 64287 \\
	\texttt{sascha.eisentraeger@tu-darmstadt.de} \\
		\And
	\href{https://orcid.org/0000-0001-8997-3818}{\includegraphics[scale=0.06]{orcid.pdf}\hspace{1mm}Daniel Juhre} \\
	Institute~of~Mechanics \\
	Otto von Guericke University Magdeburg \\
	Magdeburg, 39106 \\
	\texttt{daniel.juhre@ovgu.de} \\
		\And
	\href{https://orcid.org/0000-0002-2162-3675}{\includegraphics[scale=0.06]{orcid.pdf}\hspace{1mm}Alexander D\"uster} \\
	Institute for Ship Structural Design and Analysis (M-10) \\
	Numerical Structural Analysis with Appl. in Ship Technology \\
	Hamburg University of Technology \\
	Hamburg, 21073 \\
	\texttt{alexander.duester@tuhh.de} \\
}

\renewcommand{\undertitle}{Preprint submitted to \textit{Journal of Sound and Vibration}}
\renewcommand{\shorttitle}{Vibroacoustics using a ficticious domain approach}

\hypersetup{
pdftitle={Vibroacoustic simulations of acoustic damping materials using a fictitious domain approach},
pdfsubject={cs.NA, cs.CE},
pdfauthor={Lars Radtke},
pdfkeywords={vibroacoustics, ficticious domain approach, monolithic coupling, explicit dynamics, acoustic damping materials},
}

\begin{document}
\maketitle

\begin{abstract}
	The numerical investigation of acoustic damping materials, such as foams, constitutes a valuable enhancement to experimental testing.
Typically, such materials are modeled in a homogenized way in order to reduce the computational effort and to circumvent the need for a computational mesh that resolves the complex micro-structure.
However, to gain detailed insight into the acoustic behavior, e.g., the transmittance of noise, such fully resolved models are mandatory.
The meshing process can still be avoided by using a ficticious domain approach.
We propose the finite cell method, which combines the ficticious domain approach with high-order finite elements and resolves the complex geometry using special quadrature rules.
In order to take into account the fluid-filled pores of a typical damping material, a coupled vibroacoustic problem needs to be solved.
To this end, we construct two separate finite cell discretizations and prescribe coupling conditions at the interface in the usual manner. 
The only difference to a classical boundary fitted approach to vibroacoustics is that the fluid-solid interface is immersed into the respective discretization and does not correspond to the element boundaries.
The proposed enhancement of the finite cell method for vibroacoustics is verified based on a comparison with commercial software and used within an exemplary application.
\end{abstract}


\section{\uppercase{Introduction}}
\label{sec:introduction}

The numerical analysis of foamed acoustic damping materials is a challenging task that requires the application of non-standard modeling and discretization methods.
In particular the complex micro-structure renders standard finite element discretization unsuitable due to an unacceptably difficult mesh generation in combination with an unacceptably large number of elements needed to resolve the fine geometric details.
Typically, such structures are therefore modeled in a homogenized way, where the effective material parameters are determined based on experiments, see e.g., \cite{spannan2020,duvigneau2016}.
In order to support experimental investigations and replace them partly by simulations, a numerical model that is based on the fully resolved geometry is of great interest.
It is further desired that the simulations are done in the time domain, such that viscoelastic effects can be directly taken into account.
To this end, this work utilizes the finite cell method (FCM), a combination of a fictitious domain method with high-order finite elements. A vibroacoustic coupling is realized by superimposing two separate FCM discretizations, one for the structure (the foam material) and one for the fluid (the foam pores).

The FCM was first proposed in \cite{parvizian2007} and \cite{duester2008_1}. Several similar discretization approaches can be found in literature and are oftentimes known under a different name.
Early references used the terms \textit{fixed grid FEM} (FGFEM) \cite{garcia1999}, \textit{unfitted FEM} \cite{hansbo2002}, or more generally \textit{immersed boundary methods} \cite{roma1999}.
More recent publications include the terms \textit{CutFEM} \cite{burman2010,burman2012}, \textit{extended or generalized FEM} (XFEM/GFEM) \cite{fries2010}, and \textit{Cartesian grid FEM} (cgFEM) \cite{nadal2013,munoz2022}.
Common to all methods is the idea to use a Cartesian grid that discretizes the bounding box of the actual computational domain.
This so called extended domain is the union of the actual computational domain (denoted as \textit{physical domain}) and a \textit{fictitious domain}.
The actual geometry is accounted for during the intergration of the underlying weak form of the problem, where integrands are multiplied by an indicator function. The indicator function in the fictitious domain is set to zero or a very small positive number, while it is set to one in the physical domain.

The FCM has been successfully applied in the context of several problem types in the field of structural mechanics. This includes small and large deformation analyses \cite{schillinger2012} and plasticity \cite{abedian2013}. Also in the scope of dynamic problems, the FCM has been applied previously using either implicit \cite{elhaddad2015} or explicit \cite{joulaian2014, duczek2014} time integration schemes.
A drawback of the presence of cut cells, i.e., cells, which are cut by the boundary of the physical domain is that instabilities may arise. 
Especially cells with little support in the physical domain are problematic, as they increase the condition number of the resulting system matrices.
This yields poor convergence or even divergence of iterative linear system solvers for static problems and reduces the critical time step size of explicit time integration schemes for dynamic problems.
To this end, several remedies have been proposed, including tailored preconditioners \cite{dePrenter2019} and eigenvalue stabilization techniques \cite{garhuom2022}.
Problems arising for large deformation analyses can be counteracted by the remeshing strategy proposed in \cite{garhuom2020}.

The present work builds upon the classical FCM without any of the mentioned remedies.
This allows for an investigation of the suitability and efficiency of the classical FCM in the context of coupled vibroacoustic problems, which -- to the best of the authors' knowledge -- has not been done before.
In this context, the foam material is modeled as a linear elastic solid and the discretization follows \cite{joulaian2014, duczek2014}, however, no mass lumping is applied.
The pores are modeled as an acoustic fluid and the discretization follows \cite{buerchner2023}.
The central difference method is used for time integration. 
The use of an explicit scheme with an inherent critical time step size is justified by the strict requirements on the time step size posed by the physics.

The remainder of this work is structured as follows. In Sec.~\ref{sec:methodology} we introduce the monolithically coupled physical model and its discretization using the FCM.
In Sec.~\ref{sec:verification} we propose a benchmark problem and compare the simulation results obtained with our implementation of the FCM to results obtained with a commercial software that uses a boundary fitted FEM-based discretization.
In Sec.~\ref{sec:application} an exemplary application is considered that highlights the need for a fully resolved model and in Sec.~\ref{sec:conclusion} conclusions are drawn and an outlook to future works is given.

\section{\uppercase{Methodology}}
\label{sec:methodology}
In this section, we first recall the basic equations governing the interaction of an elastic solid and an acoustic fluid.
Afterwards, we introduce the discretization using the finite cell method.

\subsection{Vibroacoustic coupling}
We consider a linear elastic structure within the domain $\Omega^\mathrm{s}$. 
Its mechanical behavior is governed by the balance of linear momentum
\begin{align}
\rho^\mathrm{s}\, \ddot{\bm{d}} - \mathrm{div}(\bm{\sigma}) = \bm{b} \label{eq:structure_pde},
\end{align}
where $\rho^\mathrm{s}$ is the mass density, $\bm{d}(\bm{x}, t)$ is the displacement field, $\bm{\sigma}$ is the Cauchy stress tensor and $\bm{b}$ is a volumetric load.
The stress tensor $\bm{sigma}$ is related to the engineering strain $\bm{\varepsilon}$ through Hooke's law, i.e.,
\begin{align}
\bm{\sigma} = \mathcal{C} \, \bm{\varepsilon}, \quad \text{with } \bm{\varepsilon} = \frac{1}{2} \left( \nabla \bm{d} + (\nabla \bm{d})^\mathrm{T} \right),
\end{align}
where $\mathcal{C}$ is the fourth-order elasticity tensor.
For a more detailed introduction to structural dynamics in the context of finite elements we refer to \cite{bathe1996,zienkiewicz2000,hughes2012}.

We further consider an acoustic fluid that occupies the domain $\Omega^\mathrm{f}$.
Its behavior is governed by the scalar wave equation
\begin{align}
\ddot{\Psi} - \nabla \cdot \left( c^2 \, \nabla \Psi \right) = f,
\label{eq:fluid_pde}
\end{align}
where $\Psi(\bm{x}, t)$ denotes the velocity potential, $c$ is the wave velocity and $f$ is a volumetric load.
The acoustic pressure and the particle velocity are given as
\begin{align}
p^\mathrm{f} = \rho^\mathrm{f} \, \dot{\Psi} \quad \text{and} \quad
\bm{v}^\mathrm{f} = - \nabla \Psi,
\end{align}
where $\rho^\mathrm{f}$ denotes the fluid density.
      
In order to establish a mechanical coupling between the two domains, the following coupling conditions are prescribed along the interface $\Gamma^\mathrm{i} = \Omega^\mathrm{s} \cap \Omega^\mathrm{f}$:
\begin{align}
-\bm{\sigma} \, \bm{n} &= p^\mathrm{s} \, \bm{n} = p^\mathrm{f} \, \bm{n},\label{eq:coupling_condtion_structure}
\\
\bm{v}^\mathrm{f} \cdot \bm{n} &= \dot{\bm{d}} \cdot \bm{n},
\label{eq:coupling_condition_fluid}
\end{align}
where $\bm{n}$ denotes the unit normal vector of $\Gamma^\mathrm{i}$ that points out of  $\Omega^\mathrm{s}$ and into $\Omega^\mathrm{f}$.

Multiplication of equations \eqref{eq:structure_pde} and \eqref{eq:fluid_pde} by a test function ($\delta\bm{d}$ and $\delta \Psi$, respectively), integration over the domains, and integration by parts of the divergence terms results in a weak form of the two problems:
\begin{align}
\int_{\Omega^\mathrm{s}} \rho^\mathrm{s}\, \ddot{\bm{d}} \cdot \delta \bm{d} \, \mathrm{d}\Omega^\mathrm{s} 
+
\int_{\Omega^\mathrm{s}} \bm{\varepsilon} \cdot \mathcal{C} \,  \delta\bm{\varepsilon} \, \mathrm{d}\Omega^\mathrm{s} 
=
\int_{\Gamma^\mathrm{s,N}} \bm{\sigma}\,\bm{n} \cdot \delta \bm{d} \, \mathrm{d}\Gamma^\mathrm{s,N}
+
\int_{\Omega^\mathrm{s}} \bm{b} \cdot \delta \bm{d} \, \mathrm{d}\Omega^\mathrm{s},
\end{align}
\begin{align}
\int_{\Omega^\mathrm{f}} \ddot{\Psi} \cdot \delta \Psi \, \mathrm{d}\Omega^\mathrm{f} 
+
\int_{\Omega^\mathrm{f}} c^2 \, \nabla \Psi \cdot \nabla \delta \Psi \, \mathrm{d}\Omega^\mathrm{f} 
=
- \int_{\Gamma^\mathrm{f,N}} c^2 \, \nabla \Psi \cdot \bm{n} \, \delta \Psi \, \mathrm{d}\Gamma^\mathrm{f,N}
+
\int_{\Omega^\mathrm{f}} f \, \delta \Psi \, \mathrm{d}\Omega^\mathrm{f}, 
\end{align}
where $\delta\bm{\varepsilon}$ denotes the variation of $\bm{\varepsilon}$.
We have introduced here $\Gamma^\mathrm{f,N}$ and $\Gamma^\mathrm{s,N}$ as those parts of $\partial \Omega^\mathrm{f}$ and $\partial \Omega^\mathrm{s}$, where Neumann boundary conditions are prescribed. 
This includes the coupling interface $\Gamma^\mathrm{i}$ but excludes parts, where Dirichlet boundary conditions are prescribed, which are enforced in a strong sense. This is only possible, because Dirichlet boundary conditions are prescribed only on element boundaries for the problems considered here.
It is noted that Dirichlet boundary conditions have to be enforced weakly, if they are prescribed on an immersed interface to avoid locking. However, this is not the case for the examples considered in this work.

\subsection{Discretization}
\begin{figure}
    \centering
    \includegraphics[width=0.8\textwidth]{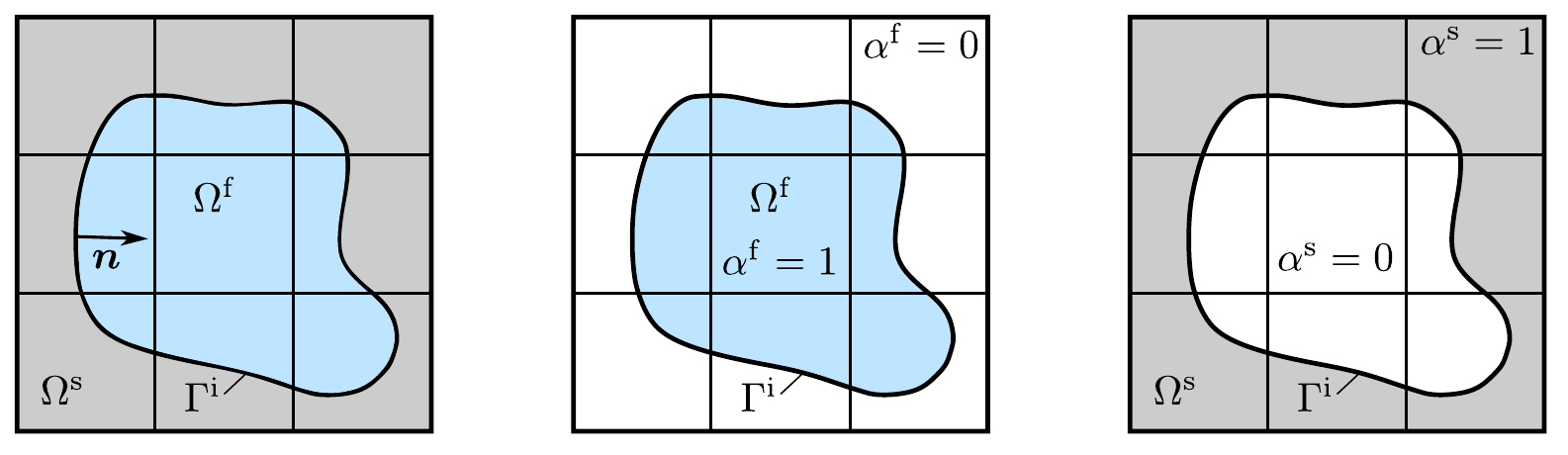}
    \caption{Discretization of the fluid domain $\Omega^\mathrm{f}$ and the structure domain $\Omega^\mathrm{s}$ using two independent FCM discretizations.}
    \label{fig:fcm_coupling}
\end{figure}
In order to arrive at a discretization in the sense of the FCM as illustrated in Fig.~\ref{fig:fcm_coupling}, we extend the integrals in the weak form and introduce the indicator functions 
\begin{align}
\alpha^\mathrm{s}(\bm{x}) = \left\{ \begin{array}{ll}
   1 & \text{ for } \bm{x} \in \Omega^\mathrm{s}, \\
   \alpha^\mathrm{min}   & \text{ else }
\end{array} \right.
\quad \text{ and } \quad
\alpha^\mathrm{f}(\bm{x}) =  \left\{ \begin{array}{ll}
   1 & \text{ for } \bm{x} \in \Omega^\mathrm{f}, \\
   \alpha^\mathrm{min}   & \text{ else }.
\end{array} \right.
\end{align}
If not stated otherwise, we set $\alpha^\mathrm{min}=10^{-8}$.
After including the coupling conditions according to Eqs.~\eqref{eq:coupling_condtion_structure}~and~\eqref{eq:coupling_condition_fluid} and assuming that $\bm{\sigma}\,\bm{n}=\bm{t}$ on $\Gamma^\mathrm{s,N} \setminus \Gamma^\mathrm{i}$ and that $\bm{v}^\mathrm{f} \cdot \bm{n} = 0$ on $\Gamma^\mathrm{f,N} \setminus \Gamma^\mathrm{i}$ we obtain
\begin{align}
\int_{\Omega^\mathrm{e}}  \alpha^\mathrm{s} \, \rho^\mathrm{s}\, \ddot{\bm{d}} \cdot \delta \bm{d} \, \mathrm{d}\Omega^\mathrm{e} 
+
\int_{\Omega^\mathrm{e}}  \alpha^\mathrm{s} \,  \bm{\varepsilon} \cdot \mathcal{C} \,  \delta\bm{\varepsilon} \, \mathrm{d}\Omega^\mathrm{e} 
+
\int_{\Gamma^\mathrm{i}} \rho^\mathrm{f} \,\dot{\Psi} \, \bm{n} \cdot \delta \bm{d} \, \mathrm{d}\Gamma^\mathrm{i}
=
\int_{\Omega^\mathrm{e}} \alpha^\mathrm{s} \, \bm{b} \cdot \delta \bm{d} \, \mathrm{d}\Omega^\mathrm{e} 
+
\int_{\Gamma^\mathrm{s}} \bm{t} \cdot \delta \bm{d} \, \mathrm{d}\Gamma^\mathrm{s,N} 
\label{eq:weak_form_structure}
\end{align}
and
\begin{align}
\int_{\Omega^\mathrm{e}} \alpha^\mathrm{f} \, \ddot{\Psi} \cdot \delta \Psi \, \mathrm{d}\Omega^\mathrm{e} 
+
\int_{\Omega^\mathrm{e}}  \alpha^\mathrm{f} \,  c^2 \, \nabla \Psi \cdot \nabla \delta \Psi \, \mathrm{d}\Omega^\mathrm{e} 
-
\int_{\Gamma^\mathrm{i}} c^2  \, \bm{n} \cdot \dot{\bm{d}} \, \delta \Psi \, \mathrm{d}\Gamma^\mathrm{i} 
=
\int_{\Omega^\mathrm{e}} \alpha^\mathrm{f} \, f \, \delta \Psi \, \mathrm{d}\Omega^\mathrm{e}.
\label{eq:weak_form_fluid}
\end{align}
Therein, $\Omega^\mathrm{e}$ denotes the extended domain. For the illustrative problem shown in Fig.~\ref{fig:fcm_coupling} and for all simulations in this work, $\Omega^\mathrm{e}$ is equivalent for fluid and structure. 
Substituting the usual approximation based on finite element shape functions
\begin{align}
    \bm{d} \approx \sum_{i=1}^{n^\mathrm{s}} N^\mathrm{s}_i \, \mathbf{U}^\mathrm{s}_i = \mathbf{N}^\mathrm{s} \, \mathbf{U}^\mathrm{s}
    \quad \text{and} \quad
    \Psi \approx \sum_{i=1}^{n^\mathrm{f}} N^\mathrm{f}_i \, U^\mathrm{f}_i = \mathbf{N}^\mathrm{f} \, \mathbf{U}^\mathrm{f}
\end{align}
into Eqs.~\eqref{eq:weak_form_structure} and \eqref{eq:weak_form_fluid} we obtain the global system of equations
\begin{align}    
\left[
\begin{array}{cc}
\mathbf{M}^\mathrm{s} & \bm{0} \\
\bm{0} & \mathbf{M}^\mathrm{f}
\end{array}
\right]
\left[
\begin{array}{c}
\ddot{\mathbf{d}} \\
\ddot{\bm{\Psi}}
\end{array}
\right]
+
\left[
\begin{array}{cc}
\bm{0} & \mathbf{C}^\mathrm{s} \\
\mathbf{C}^\mathrm{f} & \bm{0}
\end{array}
\right]
\left[
\begin{array}{c}
\dot{\mathbf{d}} \\
\dot{\bm{\Psi}}
\end{array}
\right]
+
\left[
\begin{array}{cc}
\mathbf{K}^\mathrm{s} & \bm{0} \\
\bm{0} & \mathbf{K}^\mathrm{f}
\end{array}
\right]
\left[
\begin{array}{c}
\mathbf{d} \\
\bm{\Psi}
\end{array}
\right]
=
\left[
\begin{array}{c}
\mathbf{f}^\mathrm{s} \\
\mathbf{f}^\mathrm{f}
\end{array}
\right].
\label{eq:global_system}
\end{align}
Therein,
\begin{align}
    \mathbf{M}^\mathrm{s} = \int_{\Omega^\mathrm{e}} \alpha^\mathrm{s} \, \rho^\mathrm{s} \, {\mathbf{N}^\mathrm{s}}^\mathrm{T} \, \mathbf{N}^\mathrm{s} \, \mathrm{d}\Omega^\mathrm{e} \, 
\quad \text{ and } \quad
    \mathbf{M}^\mathrm{f} = \int_{\Omega^\mathrm{e}} \alpha^\mathrm{f} \, {\mathbf{N}^\mathrm{f}}^\mathrm{T} \, \mathbf{N}^\mathrm{f} \, \mathrm{d}\Omega^\mathrm{e} \, 
\end{align}
denote the mass matrices.
For the  coupling matrices we obtain  $\mathbf{C}^\mathrm{s} = \rho^\mathrm{f} \, \mathbf{C}$ and $\mathbf{C}^\mathrm{f} = - c^2 \, \mathbf{C}^\mathrm{T}$, where
\begin{align}
    \mathbf{C} = \int_{\Gamma^\mathrm{i}} {\mathbf{N}^\mathrm{s}}^\mathrm{T} \, \bm{n} \, \mathbf{N}^\mathrm{f} \, \mathrm{d}\Gamma^\mathrm{i}.
\end{align}
The stiffness matrices are given as
\begin{align}
    \mathbf{K}^\mathrm{s} = \int_{\Omega^\mathrm{e}} \alpha^\mathrm{s} \, {\mathbf{B}}^\mathrm{T} \, 
    \mathcal{C}^\mathrm{V} \, \mathbf{B} \, \mathrm{d}\Omega^\mathrm{e} \, 
\quad \text{ and } \quad
    \mathbf{K}^\mathrm{f} = \int_{\Omega^\mathrm{e}} \alpha^\mathrm{f} \, c^2 \, \mathbf{G}^\mathrm{T} \, \mathbf{G} \, \mathrm{d}\Omega^\mathrm{e}, 
\end{align}
where the strain-displacement matrix $\mathbf{B}$ and the gradient-potential matrix $\mathbf{G}$ are defined such that
\begin{align}
    \bm{\varepsilon}^\mathrm{V} \approx \mathbf{B} \, \mathbf{U}^\mathrm{s}
\quad \text{ and } \quad
    \nabla \Psi \approx \mathbf{G} \, \mathbf{U}^\mathrm{f}.
\end{align}
In the above definitions the matrix $\mathcal{C}^\mathrm{V}$ and the vector $\bm{\epsilon}^\mathrm{V}$ denote the elasticity tensor and the strain tensor in Voigt notation, which yields
\begin{align}
{\delta \mathbf{U}^\mathrm{s}}^\mathrm{T} \, \mathbf{B}^\mathrm{T} \mathcal{C}^\mathrm{V} \, \mathbf{B} \, \mathbf{U}^\mathrm{s}
\approx
\bm{\varepsilon}^\mathrm{V} \cdot \mathcal{C}^\mathrm{V} \delta \bm{\varepsilon}^\mathrm{V}
=
\bm{\varepsilon} \cdot \mathcal{C} \, \delta \bm{\varepsilon}. 
\end{align}
The load vectors in Eq.~\eqref{eq:global_system} are given as
\begin{align}
    \mathbf{f}^\mathrm{s} = \int_{\Omega^\mathrm{e}} \alpha^\mathrm{s} \, {\mathbf{N}^\mathrm{s}}^\mathrm{T} \, 
    \bm{b} \, \mathrm{d}\Omega^\mathrm{e} 
    +
    \int_{\Gamma^\mathrm{s}} {\mathbf{N}^\mathrm{s}}^\mathrm{T} \, 
    \bm{t} \, \mathrm{d}\Omega^\mathrm{e} 
\quad \text{ and } \quad
    \mathbf{f}^\mathrm{f} = \int_{\Omega^\mathrm{e}} \alpha^\mathrm{f} \, {\mathbf{N}^\mathrm{f}}^\mathrm{T} \, 
    f \, \mathrm{d}\Omega^\mathrm{e}.
\end{align}

If not stated otherwise, we choose $N^\mathrm{s}_i = N^\mathrm{f}_i$ to be hierarchical shape functions of degree $p$ as detailed, e.g., in \cite{szabo1991,duester2001}.
Focusing on two-dimensional benchmark problems for the time being, efficiency considerations are not of utmost importance. 
We therefore apply no lumping to the mass matrix, neither by nodal quadrature, nor by explicit lumping schemes such as row-summing or diagonal scaling (commonly referred to as HRZ lumping).
However, building upon the findings in \cite{joulaian2014,duczek2014}, we finally aim at a discretization with Lagrange polynomials in combination with a nodal quadrature.
This combination in the immersed setting was introduced as the spectral cell method (SCM) in \cite{joulaian2014,duczek2014}, referencing the spectral element method (SEM), which has been known previously for its improved performance for explicit dynamic simulations. For early works on the SEM, see \cite{patera1984,ronquist1987}, for the variant considered in the references above, see \cite{cohen2002,komatitsch2002}.
However, the choice for a lumping scheme applied to cut cells must be taken with care and in consideration of the applied quadrature rule.

At this point, we want to stress that the issue of lumping is still an open research question and until now, no fully satisfactory solution has been found for immersed methods, high-order finite elements, and isogeometric analysis (IGA). 
However, an increased interest in mass lumping is observed lately, especially in the context of IGA, which could also benefit developments in the framework of the FCM.
For recent advances in dual mass lumping, the reader is referred to \cite{held2024, nguyen2024}, while the underlying mathematical theory of mass lumping is discussed in \cite{voet2023}. It is hoped that these results will also lead to improved lumping techniques for other numerical methods, where highly accurate diagonal mass matrices have not been obtained yet.

\subsubsection{Quadrature} 
The integrals arising in the weak forms \eqref{eq:weak_form_structure} and \eqref{eq:weak_form_fluid} can no longer be accurately integrated using standard Gaussian quadrature rules.
Instead, the discontinuities introduced by the indicator functions $\alpha^\mathrm{s}$ and $\alpha^\mathrm{f}$ have to be resolved.
We use here the well established approach of a space tree partitioning (a quadtree in two dimensions and an octree in three dimensions) for all cells that are intersected by an immersed boundary.
It is noted that alternative approaches that achieve the same accuracy with fewer quadrature points have been developed.
They include smart space trees \cite{kudela2015,kudela2016}, merging of sub-cells \cite{petoe2020,petoe2022} and Boolean operations \cite{abedian2017,petoe2023}.
However, for the linear problem at hand, the computation of the system matrices is not a bottleneck, because it has to be performed only once, as opposed to nonlinear problems.
Alternative quadrature rules such as moment fitting (see \cite{hubrich2017,hubrich2019,garhuom2022_2}) are interesting in particular because they may allow to construct nodal quadrature rules for cut cells as demonstrated in \cite{nicoli2022}.
A similar nodal quadrature idea has also been proposed in the context of high-order Hermite shape functions \cite{kapuria2021}, which, however, shows the same shortcomings as seen for fictitious domain methods in that only sub-optimal rates of convergence are attainable. 
This has been conclusively discussed in \cite{eisentraeger2023} and therefore, such lumping techniques still lack accuracy and need to be further improved.
In this paper we first establish vibroacoustic simulations  using the finite cell method in combination with a consistent mass matrix and postpone the use of advanced quadrature rules that yield diagonal mass matrices by construction to future works.

\subsubsection{Geometry description}
The geometry of the problems considered in this work are defined analytically based on a levelset function $\varphi(\bm{x})$, where $\varphi<0$ in $\Omega^\mathrm{f}$ and $\varphi>0$ in $\Omega^\mathrm{s}$.
The levelset functions provide an inside-outside test, which is used to construct quadrature cells in a fully automatic way.
This paves the way for fully resolved simulations of acoustic damping materials, where the inside-outside test will be realized based on a voxel model of the geometry obtained using computed tomography (CT) scans. 
The subsequent construction of the discretization remains the same regardless of the origin of the inside-outside test and consists of the following two steps.
\begin{enumerate}
    \item The quadtree approach is used to construct quadrature cells (quadrilaterals) for the domains $\Omega^\mathrm{s}$ and $\Omega^\mathrm{f}$.
    \item A marching squares algorithm is used to construct quadrature cells (lines) for the interface $\Gamma^i$. The element boundaries are resolved by the cells. 
\end{enumerate}
For three-dimensional simulations, an octree approach is used in Step 1 and a marching cubes algorithm in Step 2.
The applicability of a spatial finite cell discretization based on CT scans of metal foams was already shown in \cite{duester2012,garhuom2022}.

\section{\uppercase{Verification benchmark}}
\label{sec:verification}
\begin{figure}
    \centering
    \includegraphics[width=\textwidth]{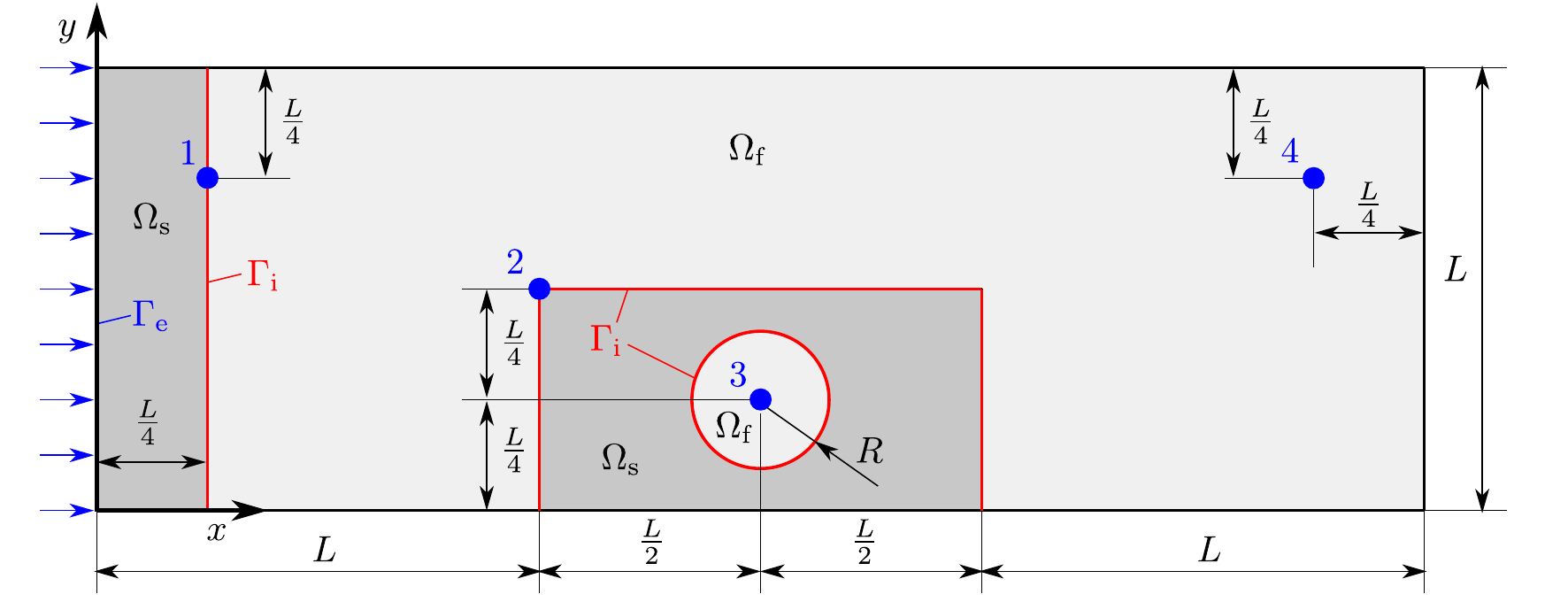}
    \caption{Geometry of the benchmark problem used for verification.}
    \label{fig:problem_sketch}
\end{figure}
In order to verify our implementation of the vibroacoustic coupling with the FCM in the in-house software \mbox{\textit{AdhoC++}} \cite{zander2014}, we consider a benchmark problem and compare the results with simulations based on the commercial FEM software \textit{Abaqus FEA}~\cite{abaqus2009}.
Figure~\ref{fig:problem_sketch} shows the geometry of the two-dimensional benchmark problem. 
We use $L=0.1\,$m and $R=\frac{L}{5}$ and assume a \textit{plane stress} condition for the structure problem.
The overall rectangular domain is decomposed into four parts, two that make up the structure domain $\Omega^\mathrm{s}$ and two that make up the fluid domain $\Omega^\mathrm{f}$.
The structure boundaries at $y=0$ and $y=L$ are clamped.
On the structure boundary at $x=0$, a pressure is prescribed, i.e., $\Gamma^\mathrm{e} = \Gamma^\mathrm{s,N}$.
All outer fluid boundaries are free (homogeneous Neumann condition).

\begin{table}
    \centering
    \caption{Constitutive parameters for structure and fluid.}
    \label{tab:properties}
    \renewcommand{\arraystretch}{1.3}
    \begin{tabular}{c|ccc|cc}
        \toprule
        & density & bulk modulus  & shear modulus & \multicolumn{2}{c}{wave velocities} \\
        &  $\rho$ (kg/m$^3$) & $\kappa$ (N/m$^2$) & $\mu$ (N/m$^2$) & pressure $c$ (m/s) & shear $c^\mathrm{s}$ (m/s) \\ \hline
        fluid & 1.225 & $0.101 \cdot 10^6$  & -- & 287.139303461 & -- \\
        structure & 50 & $\frac{5}{6} \cdot 10^6$ & $\frac{5}{13} \cdot 10^6$ & 164.082530828 & 87.7058019307 \\ \bottomrule
    \end{tabular}
\end{table}
The constitutive parameters along with the resulting wave velocities are given in Tab.~\ref{tab:properties}.
The fluid parameters are set according to those available for air.
The structure parameters result from a given Young's modulus of \mbox{$E=1\,$MPa} and a Poisson ratio of \mbox{$\nu = 0.3$}.
This corresponds to the range of values measured for acoustic foams in \cite{petru2017} in a homogenized manner. Of course, our method aims at fully resolved simulations of foams, where not the homogenized parameters but the actual material parameters for the solid phase of the foam will be used.
However, for this benchmark problem we do not consider a fully resolved foam in order to avoid thin-walled regions between the pores that would go along with such a structure.

\begin{figure}
    \centering
    \includegraphics[width=\textwidth]{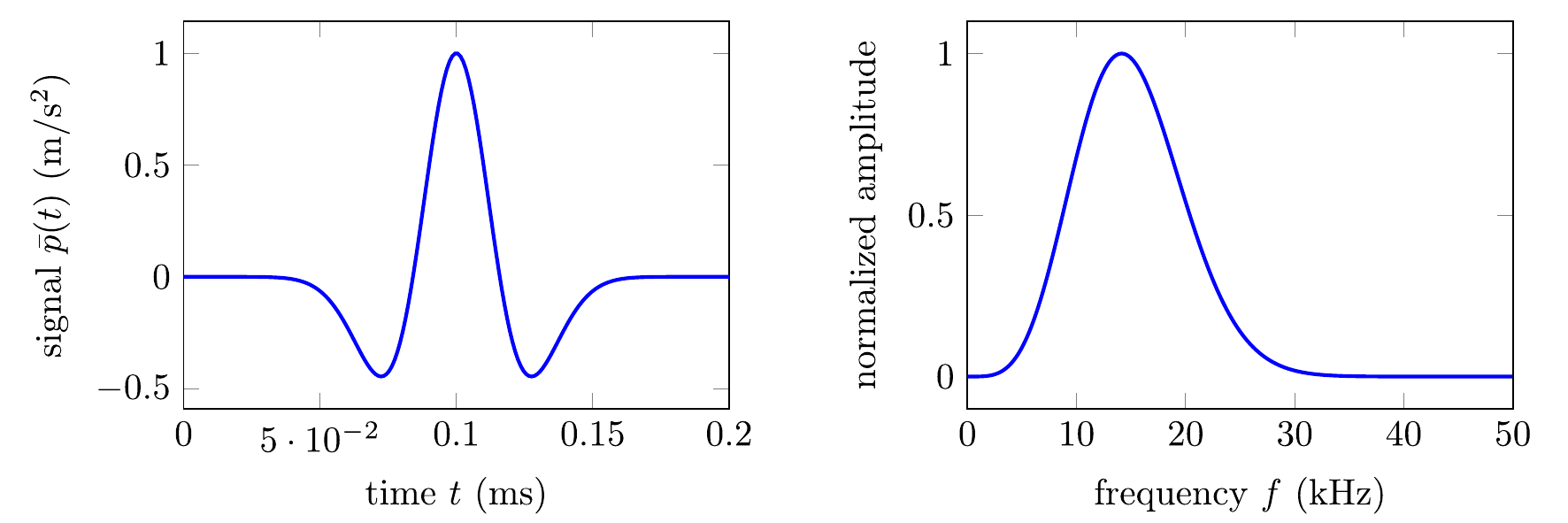}
    \caption{Excitation signal  according to Eq.~\eqref{eq:ricker} with a base frequency of $f=10^4\,$s$^{-1}$.}
    \label{fig:signals_ricker}
\end{figure}
A time varying traction $\bm{t} = \bar{p}(t) \, \bm{n}$ -- see Eq.~\eqref{eq:weak_form_structure} -- is prescribed on the left end of the structure (see Fig.~\ref{fig:problem_sketch}).
We make use of a so-called Ricker wavelet, which corresponds to a scaled second derivative of a Gaussian bell curve 
\begin{align}
\bar{p}(t) = \left( 1 - \left( \frac{t-t_0}{\sigma} \right)^2 \right) e^{\left(-\frac{(t-t_0)^2}{2 \, \sigma^2} \right)} 
\quad \text{with} \quad
\sigma= \frac{t_0}{2\,\pi}.
\label{eq:ricker}
\end{align}
The wavelet is shown in Fig.~\ref{fig:signals_ricker} for the chosen characteristic time $t=10^{-4}\,$s.
We consider a time interval with a duration of $T=2\cdot 10^{-3}$s and set the time step size to $\Delta t = 10^{-7}$s.
For comparison with the reference solution, we record the solution at the observer locations indicated in Fig.~\ref{fig:problem_sketch}.

\subsection{Discretization}
\begin{figure}
    \centering
    \includegraphics[width=\textwidth]{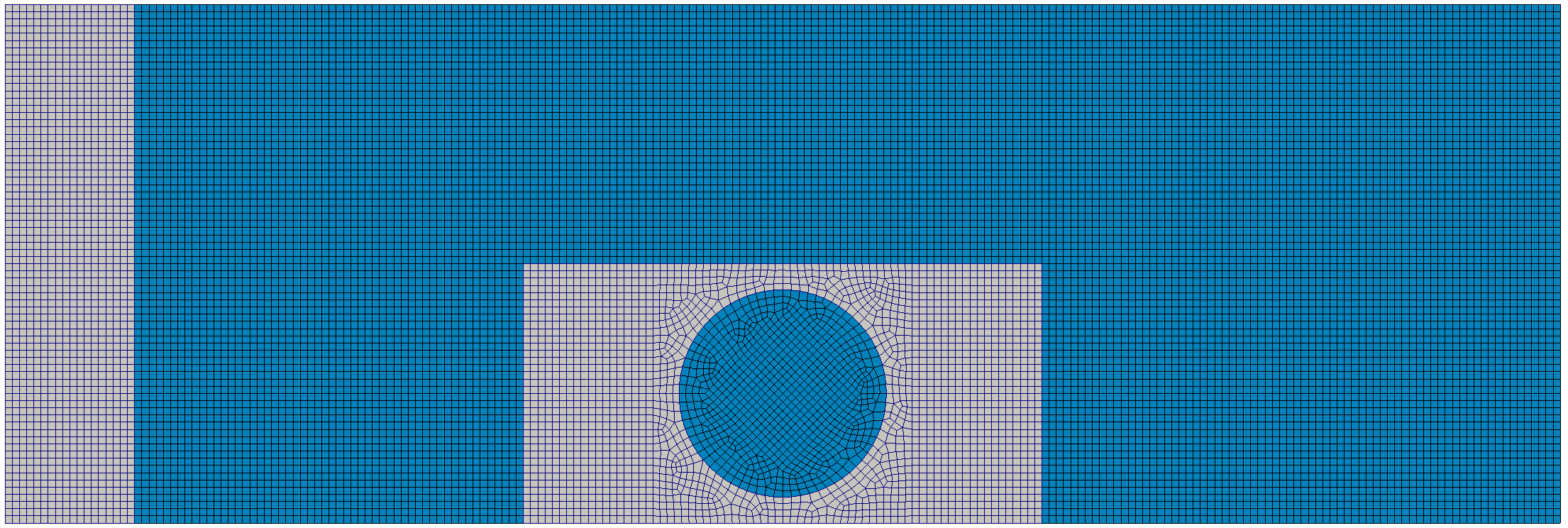}

    \vspace{1em}
    
    \includegraphics[width=\textwidth]{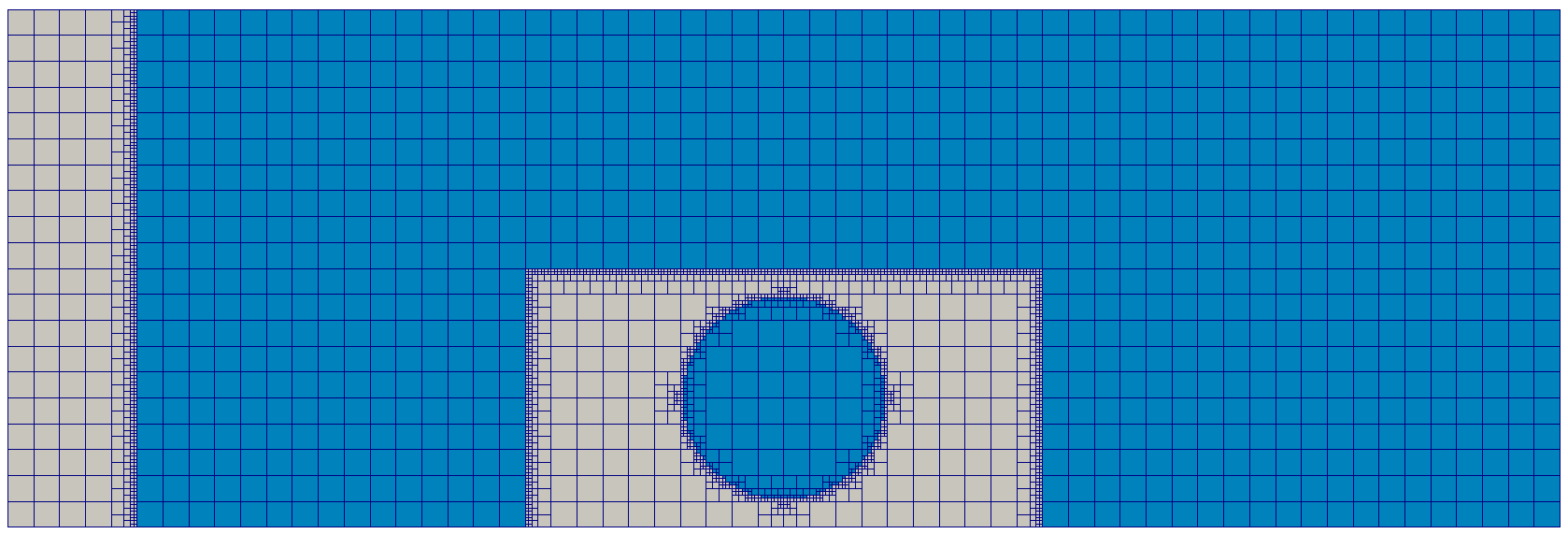}
    \caption{Top: Illustration of the boundary fitted mesh used for the reference computations. Bottom: Illustration of the FCM discretization including quadrature cells. The actually used discretizations are finer than the ones shown.}
    \label{fig:discretizations}
\end{figure}
For the reference solution computed with \textit{Abaqus FEA} \cite{abaqus2009} the domain is discretized using $34670$ quadrilateral elements for the fluid subproblem and $8985$ quadrilateral elements for the structure subproblem. 
The ansatz order is set to $p=2$, which yields $105126$ and $55472$ degrees of freedom for the fluid and the structure subproblem, respectively.
The elements are of serendipity type, i.e., eight-noded quadrilateral elements are used.
Figure~\ref{fig:discretizations} (top) shows a coarse version of the mesh for illustration.
The implicit {Hilber-Hughes-Taylor} time integration scheme {with the parameters $\alpha=0$, $\beta=0.25$ and $\gamma=0.5$ and a fixed time step size of $\Delta t = 10^{-8}s$} is used. 
It is noted that by setting $\alpha=0$ we essentially use the trapezoidal rule from the Newmark family of time integration schemes.

The FCM discretization consists of 180 by 60 quadrilateral cells with $p=3$ (for $\Omega^\mathrm{s}$ and $\Omega^\mathrm{f}$ each), which yields $109442$ and $54721$ degrees of freedom, respectively.
We use hierarchical shape functions from the so called \textit{trunk space} as detailed in \cite{szabo1991} to improve the comparability to \textit{Abaqus FEA}.
The quadrature is done using $(p+1)^2=16$ Gauss-Legendre points on each subcell. The subcells are created based on a quadtree partitioning with a depth of $p+4=7$. 
In order to evaluate the integrals over the coupling interface $\Gamma^\mathrm{i}$, the domain boundaries are discretized by linear line segments. 
These are created in a fully automatic manner using a marching cubes algorithm with resolution $10$ by $10$ for each cell \cite{lorensen1987}.

\subsection{Results}
\begin{figure}
    \centering
    \includegraphics[width=\textwidth]{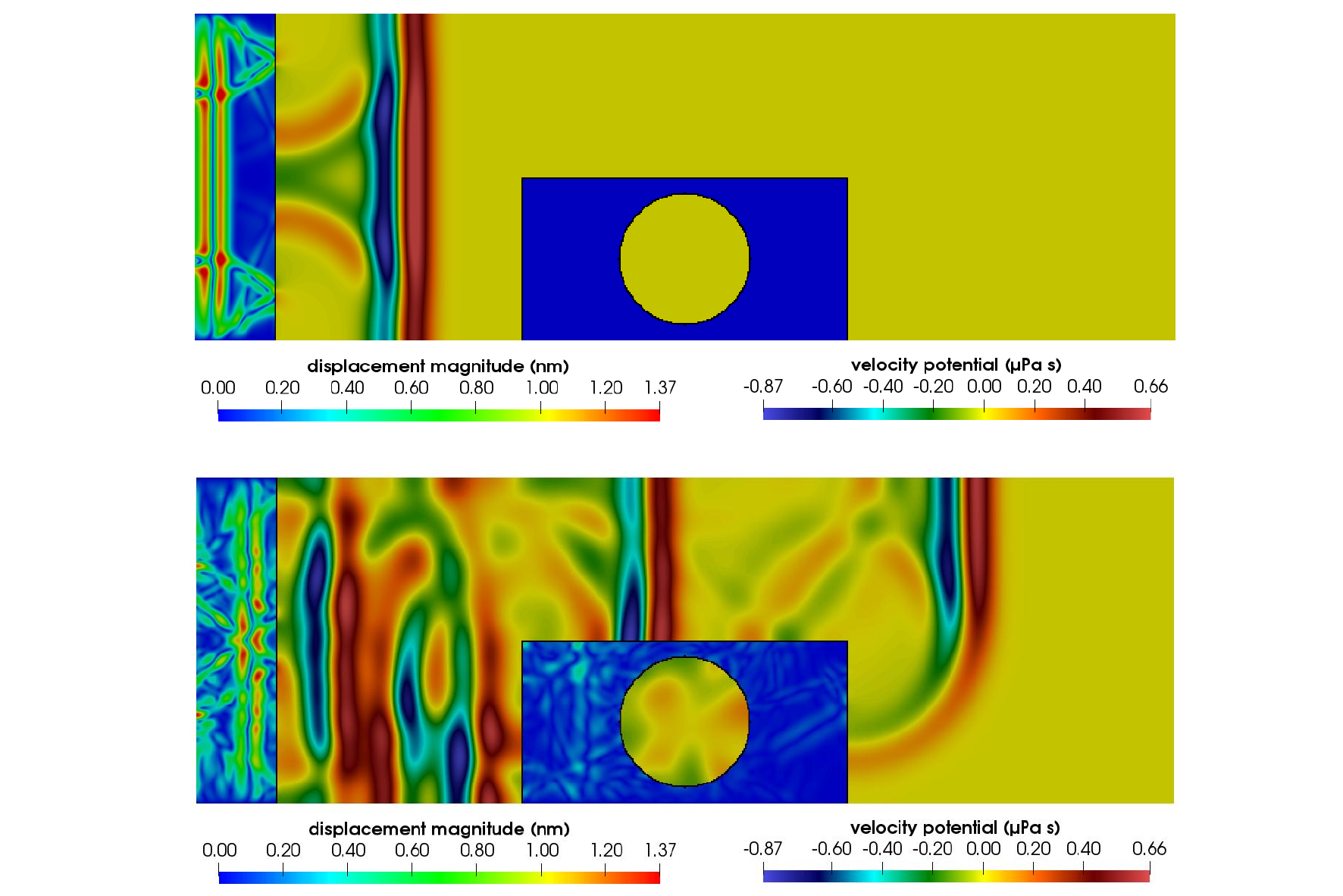}
    \caption{Overview of the simulation result. Top: solution at $t=0.4\,$ms. Bottom: solution at $t=1\,$ms.}
    \label{fig:benchmark_overview}
\end{figure}
\begin{figure}
    \centering
    \includegraphics[width=\textwidth]{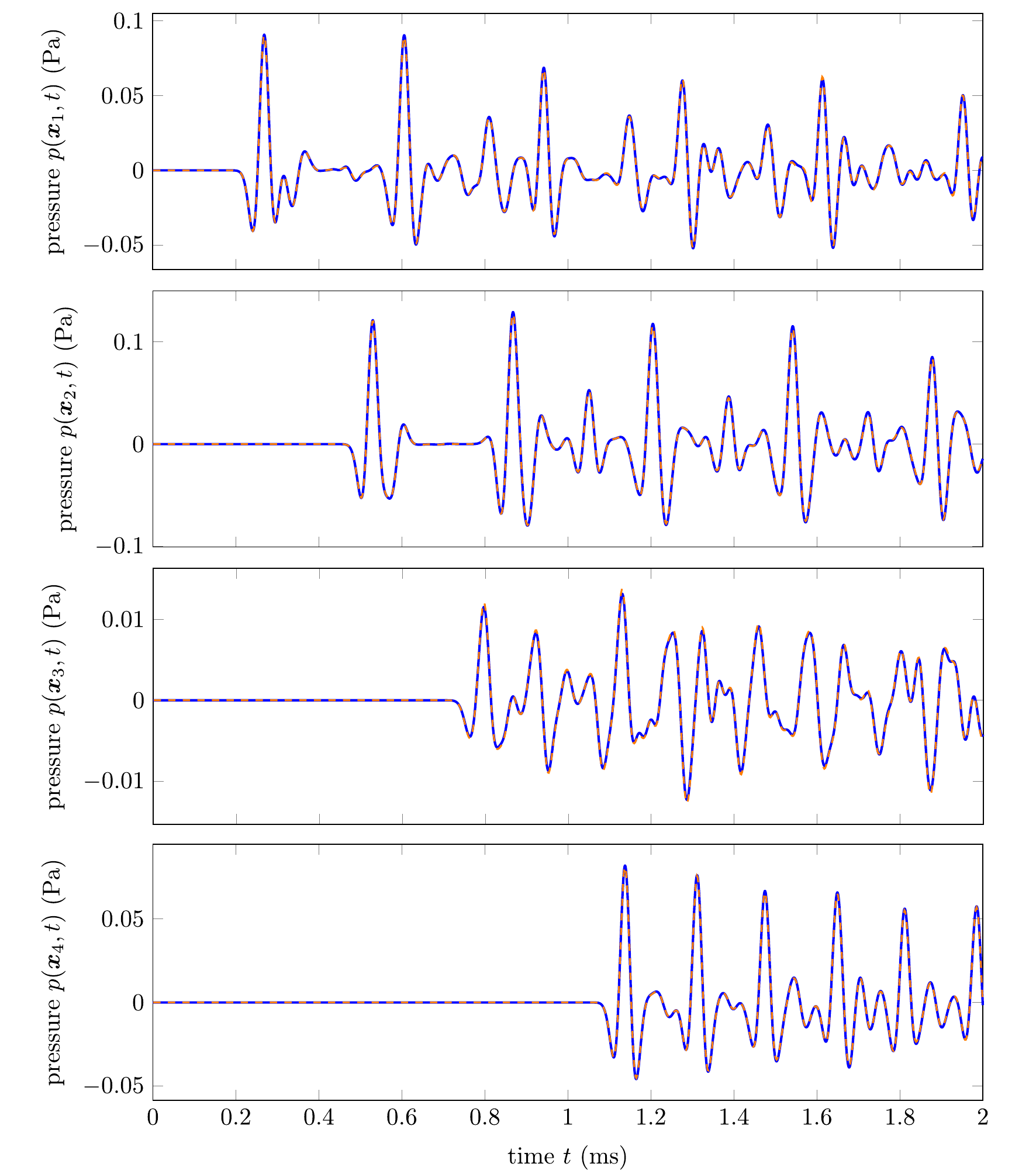}
    \caption{Comparison of the pressure signals at the observers 1--4 (see Fig.~\ref{fig:problem_sketch}). Solid lines refer to the FCM solution, dashed lines refer to the reference solution.}
    \label{fig:pressure_result}
\end{figure}
Figure~\ref{fig:benchmark_overview} shows the pressure in the fluid domain and the displacement in the structural domain at two selected time instance to give an overview of the results.
Figure~\ref{fig:pressure_result} shows the pressure recorded at the observer locations.
A very good agreement is found between the FCM results and the reference solution.
Accordingly, we view our implementation of the FCM for vibroacoustics within the framework AdhoC++ \cite{zander2014} as verified.


\section{\uppercase{Exemplary application}}
\label{sec:application}
\begin{figure}
    \centering
\includegraphics[width=\textwidth]{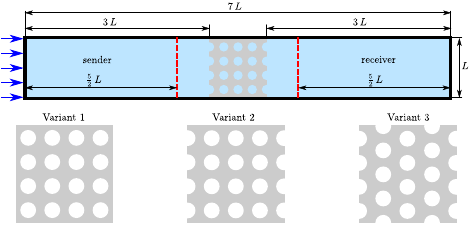}
    \caption{Top: schematic model of an impedance tube used to characterize the foam-like structure in the center region. The red lines indicate the locations, where the pressure is recorded. Bottom: three variants of the foam-like structure.}
    \label{fig:impedance_pipe}
\end{figure}

As an exemplary application we consider a two-dimensional model of an impedance tube as illustrated in Fig.~\ref{fig:impedance_pipe}.
A foam-like structure is placed in the center region and the left boundary is excited using a time varying but spatially constant Neumann boundary condition. The height of the tube is set to $L=0.05\,$m.
This yields a plane wave that travels in the positive $x$-direction and is then partly reflected at the structure.
At the indicated locations on the sender side and the receiver side, the sound pressure is recorded in order to quantify the reflection and the transmission of the structure.
The simulation is performed for three different geometries, however, all of which have the same porosity.
In a typical homogenized model, where the structure is only characterized by its porosity, the same behavior would be predicted for all geometries.
On the other hand, the fully resolved approach using the FCM is capable of identifying differences between them.

\subsection{Boundary conditions and evaluation}
For the time varying Neumann boundary condition on the left boundary, the same wavelet is used, that was also used for the benchmark case (see Fig.~\ref{fig:signals_ricker}).
This time, however, the condition is prescribed on a fluid domain boundary, resulting in a given fluid velocity, instead of a traction.
On all other fluid boundaries (except for the coupling interfaces), homogeneous Neumann boundary conditions are prescribed, which corresponds to a rigid, fully reflecting wall.
On all structure boundaries (except for the coupling interfaces) homogeneous Dirichlet boundary conditions are prescribed.

An array of $5$ observer locations is distributed equidistantly along the red dashed lines on the sender and the receiver side indicated in Fig.~\ref{fig:impedance_pipe}.
At these locations, the sound pressure is evaluated for the entire simulation time. In order to quantify the reflected and the transmitted sound, the following measures are computed.
\begin{align}
P_\mathrm{ref} = \sqrt{ \sum_j \sum_i p(\bm{x}^\mathrm{s}_i, t_j)^2 }, \quad
P_\mathrm{tra} = \sqrt{ \sum_j \sum_i p(\bm{x}^\mathrm{r}_i, t_j)^2 }.
\label{eq:rt_measures}
\end{align}
Therein, $\bm{x}^\mathrm{s}_i$ denotes the location of the $i$th observer on the sender side, $\bm{x}^\mathrm{r}_i$ denotes the location of the $i$th observer on the receiver side, and $t_j = j \, \Delta t$.

\subsection{Discretization}
\begin{figure}
    \centering
\includegraphics[width=\textwidth]{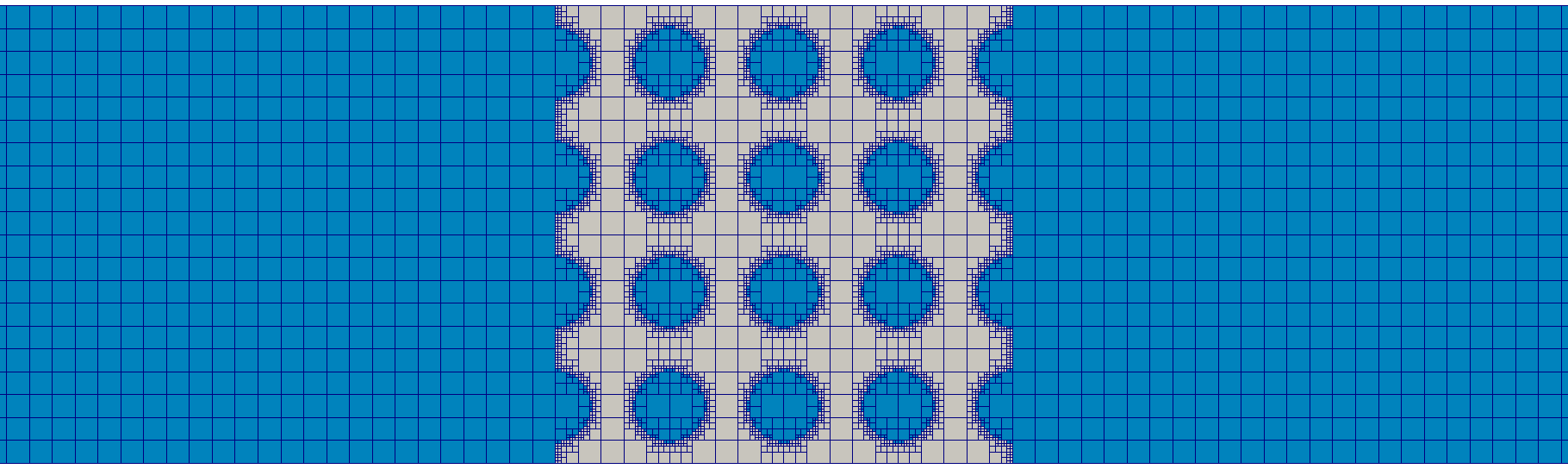}
    \caption{Quadrature cells of the FCM discretization for the problem sketched in Fig.~\ref{fig:impedance_pipe}. The actually used discretization is finer than the one shown.}
    \label{fig:discretization_impedance_pipe}
\end{figure}
The fluid domain is discretized using $60$ cells in the $y$-direction and $420$ cells in the $x$-direction, which yields $25200$ square cells in total. 
The trunk space with a polynomial degree of $p=2$ is used, which yields $76561$ degrees of freedom.
The structure is discretized using $60$ cells of the same type in $x$- and $y$-direction, which yields $3600$ cells in total and $22082$ degrees of freedom. 
For illustrative purposes, Fig.~\ref{fig:discretization_impedance_pipe} shows a zoom towards the center region of a coarser version of the FCM discretization.
A comparison with a finer discretization was performed in order to verify that the results are accurate enough.
Like for the benchmark example, a time step size of $\Delta t = 10^{-8}\,$s is used and a time interval of $2 \cdot 10^{-3}\,$s is considered in the simulation. 
For both subproblems, the quadrature is performed using $(p+1)^2=9$ Gauss-Legendre points in each direction per subcell. 
The subcells are constructed using a quadtree with a depth of $p+4=6$.

\subsection{Results}
\begin{figure}
    \centering
    \includegraphics[width=0.8\textwidth]{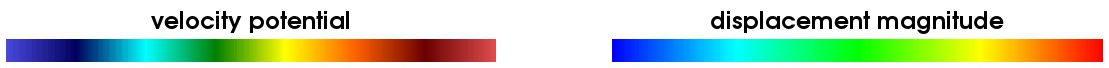}
    \vspace{1em}
    
    \includegraphics[width=0.48\textwidth]{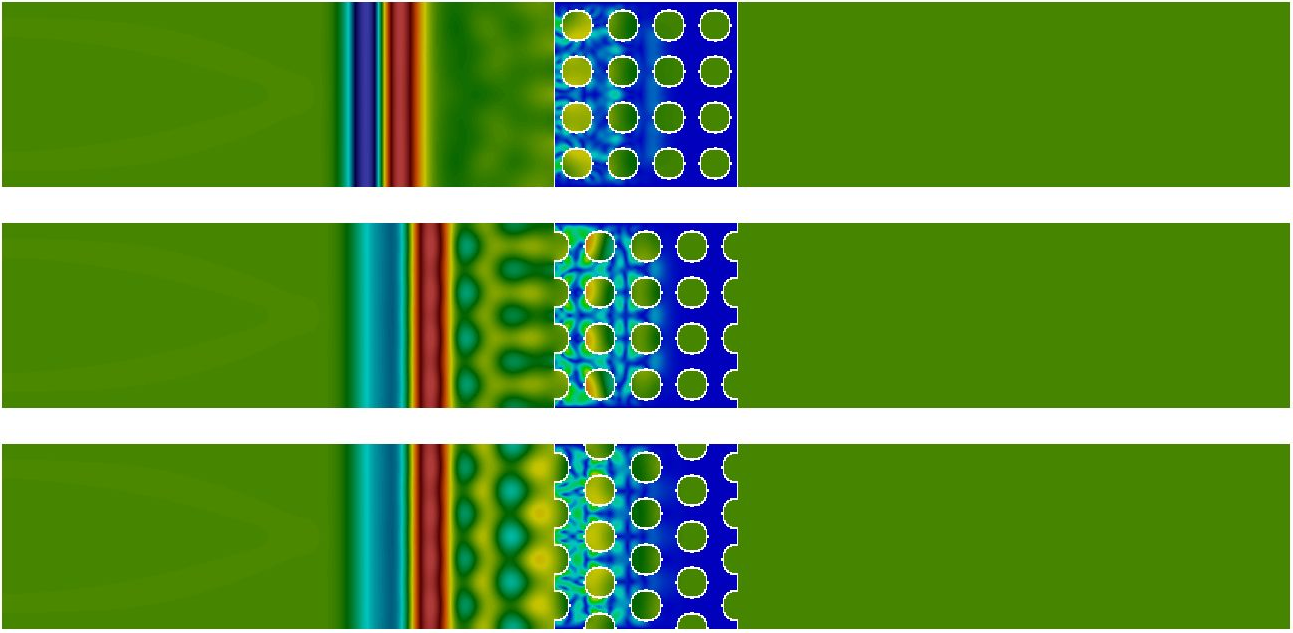}
    \includegraphics[width=0.48\textwidth]{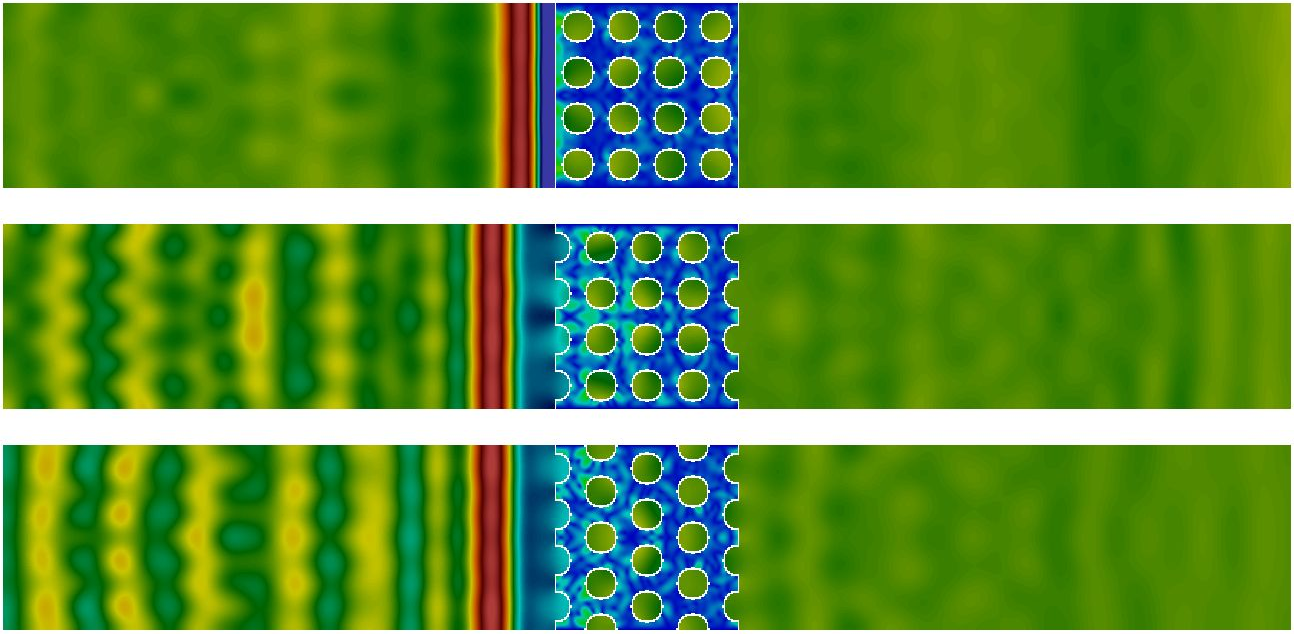}
    \caption{Snapshots of the simulation results (velocity potential in the fluid domain and displacement magnitude in the structure domain) at $t \approx 0.435\,$ms (left) and at $t \approx 1.65\,$ms (right).}
    \label{fig:fictive_foam_snapshots}
\end{figure}

In order to provide an overview of the simulation, Fig.~\ref{fig:fictive_foam_snapshots} shows a snapshot of the results at two selected time instances. At $t\approx0.435\,$ms the wave package has partly been reflected at the foam and is traveling back in the negative $x$-direction. At $t\approx 1.65\,$ms the wave front reaches the foam a second time after a reflection at the left boundary.
The influence of the foam geometry is clearly visible. While for Variant 1, the reflection is almost undisturbed, Variants 2 and 3 yield more complex patterns at the sender side.
\begin{table}
    \centering
    \caption{Expected events encountered during the simulation of the fictive foam in the impedance tube.}
    \label{tab:events}
    \renewcommand{\arraystretch}{1.3}
    \begin{tabular}{cc|l}
        \toprule
        time (ms) & distance & event  \\ \hline
        0.435 & $2.5\,L$ & A: wave front arrives at left observer locations \\
        0.609 & $3.5\,L$ & B: reflected wave front arrives at left observer locations \\
        0.784 & $4.5\,L$ & C: transmitted wave front arrives at right observer locations \\
        1.48 & $8.5\,L$ & D: twice reflected wave front arrives at left observer locations \\ 
        1.65 & $9.5\,L$ & E: three times reflected wave front arrives at left observer locations \\  
        1.65 & $9.5\,L$ & F: reflected transmitted wave front arrives at right observer locations \\ \bottomrule
    \end{tabular}
\end{table}
Table~\ref{tab:events} provides a list of the events that occur on both sides.
The given time values are computed based on the distance travelled by the wave front until reaching the location of the event.
All computations are based on the speed of sound in the fluid and therefore correspond to the earliest possible time at which the respective event can occur (the speed of sound in the structure is lower than that in the fluid, see Tab.~\ref{tab:properties}).
Accordingly, it is possible to quantify the reflectance and the transmittance of the foam variants by considering the measures introduced in Eq.~\eqref{eq:rt_measures} for particular time intervals.

\begin{figure}
    \centering
    \includegraphics[width=\textwidth]{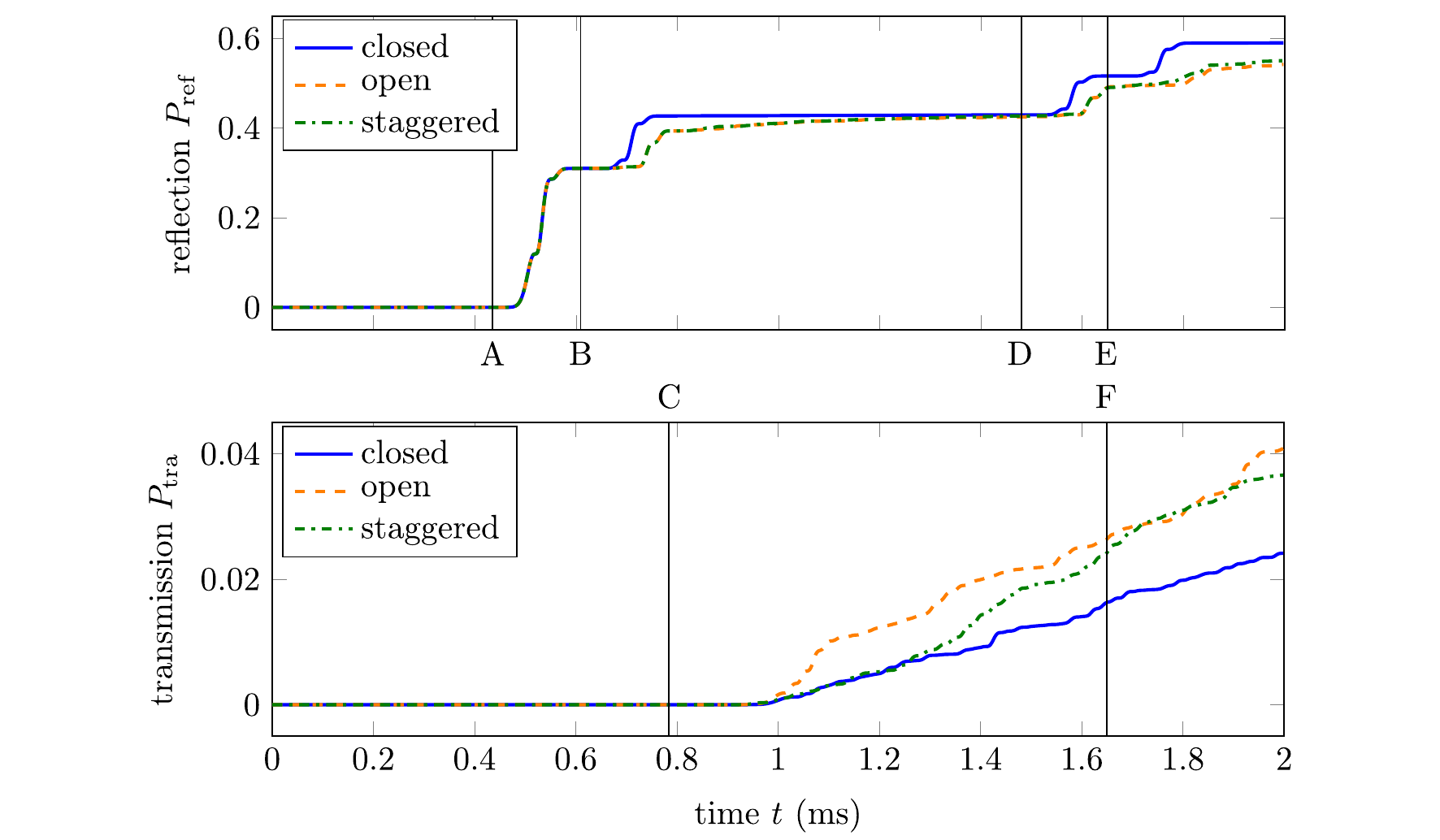}
    \caption{Reflectance and transmission measures according to Eq.~\eqref{eq:rt_measures} at the locations indicated in Fig.~\ref{fig:impedance_pipe}.}
    \label{fig:rt_measures}
\end{figure}
To this end, Fig.~\ref{fig:rt_measures} shows the two measures $P_\mathrm{ref}$ (reflectance) and $P_\mathrm{tra}$ (transmittance) over the course of the simulation.
The first increase in $P_\mathrm{ref}$ corresponds to event A, where the wave front is about to reach the observers on the left side for the first time. This increase is therefore not a result of any reflection. Instead, the second increase in $P_\mathrm{ref}$ that corresponds to event B is the result of the reflection at the foam. Any sound pressure observed before event D, when the reflected wave package arrives at the left observer locations for the third time can directly be attributed to the reflective behavior of the foam. 
Subsequent increases (events D and E) should not be taken into account. 
Considering the time interval of interest (between events B and D) a major difference between Variant 1 and the other geometries can be observed in Fig.~\ref{fig:rt_measures} (top).
Regarding the transmittance measure $P_\mathrm{trans}$ in Fig.~\ref{fig:rt_measures} (bottom), the first increase is observed shortly after event C, which corresponds to the earliest possible time, where the wave front can reach the right observer locations.
Any increase until event F (the earliest time, where the wave front can reach the right observer locations a second time after being reflected at the right boundary) can directly be attributed to the transmittance of the foam.
In accordance with the reflectance, a major difference is seen for the transmittance between Variant 1 and the other geometries. Due to the smaller scale, differences are also noticeable between Variants 2 and 3.

The results clearly show that a homogenized model purely based on the porosity is not suitable to characterize the acoustic behavior of the selected foams.
In particular, the surface geometry is found to have a severe influence on the results. 
It is noted that more realistic foams with a finer micro-structure compared to the size of the probe in an impedance tube may yield smaller differences and render a homogenized approach applicable. 
This exemplary application can therefore not be used as an argument to question such simplified modeling approaches in general, but shall be seen as a proof of concept for the developed fully resolved approach using the finite cell method.

\section{\uppercase{Summary and conclusion}}
\label{sec:conclusion}
We have developed an immersed boundary method for vibroacoustic simulations of foam-like structures in the time domain.
The implementation was verified by a comparison with a commercial software.
As an exemplary application, a simulation of an impedance tube was performed, showing the need for a fully resolved approach.
In order to enable the full potential of explicit time integration schemes and apply the method to real foams that are used as acoustic damping material, several challenges remain.

To begin with, a suitable lumping scheme has to be developed that yields a diagonal mass matrix without loss of accuracy when applied to cut cells.
While previous studies exist that present first ideas about this topic, a generally applicable approach was not yet found and a consistent mass matrix was used in the present study. Further, stabilization methods need to be developed that allow for cells with arbitrarily small support, which cannot be avoided for complex foam geometries obtained from CT scans. Without stabilization, the small support leads to large eigenvalues that restrict the time step size of conditionally stable explicit time integration schemes to an unfeasible range.

A second aspect, even though less open from a conceptual point of view, is the realization of vibroacoustic simulations with the FCM in three dimensions. In the present study, only two-dimensional models were considered and the focus was put on the coupling between the elastic structure and the acoustic fluid in an immersed setting. 

Finally, we plan to validate the simulation approach based on measurements. To this end, experiments using an impedance tube and an anechoic chamber will be conducted. The aim is to arrive at a validated simulation pipeline that takes as an input CT scans of a foam material to be investigated and yields the desired results (the acoustic damping behavior of the material) in a fully automatic way.

\section*{Acknowledgment}
The authors gratefully acknowledge the support of the DFG (Deutsche Forschungsgemeinschaft) under \mbox{DU 405/20-1} and
 \mbox{DU 1904/5-1} (grant number 503865803) and \mbox{EI 1188/3-1} (grant number 497531141).

\bibliographystyle{abbrv}
\bibliography{literature}

\end{document}